\documentclass[titlepagea4paper,11pt]{article}
\usepackage[arrow,matrix]{xy}
\usepackage{amssymb}
\usepackage{amsbsy,amsmath}
\newtheorem{thm}{Theorem}[section]
\newtheorem{lemma}[thm]{Lemma}

\newtheorem{rem}[thm]{Remark}

\newcommand{\seq}{\longrightarrow}

\newcommand{\imp}{\Rightarrow}

\newcommand{\bb}{\mathbb}

\newenvironment{proof}{\vspace{3pt}\indent
                       \textsc{\bf Proof:}\quad }
                       {\hfill$\square$\vspace{3pt}}
\newtheorem{example}[thm]{Example}
\newtheorem{construction}[thm]{Construction}
\newenvironment{ex}[1]{\begin{example}\rm{\bf #1}}{\hfill$\bigtriangleup$\end{example}}

\newenvironment{citethm}[1]{\vspace{0.3cm}\noindent{\bf Theorem \ref{#1}.} \it}{\vspace{0.3cm}}

\title{Calabi-Yau-threefolds with Picard number $\rho(X)=2$ and their K\"ahler cone I}
\author{Marco K\"uhnel\footnote{The author acknowledges gratefully support by the DFG priority program 'Global Methods in Complex Geometry'.}\\  Department of Mathematics\\  University of Bayreuth\\ D-95440 Bayreuth}
\date{\today}

\begin{document}

\maketitle

\begin{abstract}We describe properties of the K\"ahler cone of general Calabi-Yau-threefolds with Picard number $\rho(X)=2$ and prove the rationality of
the K\"ahler cone, if $X$ is a Calabi-Yau-hypersurface in a ${\bb P}^2$-bundle over ${\bb P}^2$ and
$c_3(X)\le -54$. Without the latter assumption we prove the positivity of $c_2(X)$.
\end{abstract}

\section{Introduction}In this paper, a Calabi-Yau-threefold is a compact complex K\"ahler manifold of dimension three with $K_X={\cal O}_X$ and
$H^1({\cal O}_X)=0$. 

Wilson stated in 1994 \cite{wi94} a conjecture about the rationality of the K\"ahler cone of a Calabi-Yau-threefold. It says that the K\"ahler cone
of a Calabi-Yau-threefold $X$ is rational and finitely generated in $N^1(X)$, if $c_2(X)$ is positive, i.e. $D.c_2(X)>0$ for every nef divisor $D$.

In the present paper we deal with the case of Picard number $\rho(X)=2$. 
If $X\seq S$ is an elliptic fibration onto a smooth surface and $\rho(X)=2$, then $S\cong{\bb P}^2$. 
Oguiso proves in \cite{og} that for every elliptic fibration $X\seq S$, the surface $S$, smooth or not, 
is rational.
Moreover, it is an easy argument to verify that
$C={\bb P}^1$, if $X\seq C$ is a fibration onto a normal curve $C$ and $X$ a Calabi-Yau manifold. So it is natural to consider projective spaces as base spaces of 
Calabi-Yau-fibrations. 

A rather complete picture about the K\"ahler cone should be obtained for those $X$ which can be embedded in ${\mathbb P}^n$-
bundles over ${\mathbb P}^m$ as hypersurfaces. Here we are considering the case $(n,m)=(2,2)$, i.e. 
$X$ is a family of cubic elliptic curves.  
We will prove the rationality of the K\"ahler cone for most cases:

\begin{citethm}{rat}Let ${\cal E}\seq{\mathbb P}^2$ be a rank-3-bundle, $Z:={\mathbb P}({\cal E})$, $X\subset|-K_Z|$ a
Calabi-Yau-manifold with $\rho(X)=2$ and $h^0(-K_Z)>1$. Then $\partial K(X)$ is rational.
\end{citethm}

There is also a topological criterion on $X$ 
for the condition $h^0(-K_Z)>1$, which in particular includes all possible $Z$
with $-K_Z$ nef (cf. Remark \ref{g>-18}):

\begin{citethm}{ratgamma}Let ${\cal E}\seq{\bb P}^2$ a rank-3-bundle, $X\subset{\bb P}({\cal E})=:Z$ a
Calabi-Yau threefold with  $\rho(X)=2$. If $c_3(X)\le -54$, then $h^0(-K_Z)>1$, in particular $\partial K(X)$ is rational.
\end{citethm}

To prove the rationality of $\partial K(X)$, we investigate, if
$$K(X)=K(Z)|X$$
holds, when $Z$ denotes the ${\bb P}^n$-bundle. A result of Koll\'ar states, that this is right, if $Z$ is Fano. We generalize this assumption and describe the situation,
when $K(X)\not=K(Z)|X$. 

\begin{citethm}{genkollar}Let $X\subset{\mathbb P}({\cal E})=:Z$, with ${\cal E}\seq{\mathbb P}^2$ being a rank-3-bundle, and $\rho(X)=2$. Then
$$K(X)=K(Z)|X$$
holds, unless: 

$-K_Z$ is big and nef, but not ample, $-K_Z|X$ is ample and there is a surface $G\subset Z$ such that $X\cap G=\emptyset$ and
$$[\mu G]=9{\cal O}_Z(1)^2-(6c_1({\cal E}).h+9){\cal O}_Z(1).p^*h+(9c_2({\cal E})+3c_1({\cal E}).h+9-2c_1^2({\cal E}))F$$
for a certain $\mu>0$.
\end{citethm}

This part uses the bundle-situation only for the details of the description.
Furthermore, we show that $\partial K(Z)$ is rational, if $h^0(-K_Z)>1$. Hence, under this assumption, it remains to discuss the cases
$K(X)\not=K(Z)|X$ to prove rationality of $\partial K(X)$. This is done by using some vector bundle theory, hence is essentially connected to the embedding into
a ${\bb P}^2$-bundle over ${\bb P}^2$. 

In the last part, we show

\begin{citethm}{c2>0}Let $X\subset{\bb P}({\cal E})$ a Calabi-Yau threefold with $\rho(X)=2$, where ${\cal E}\seq{\bb P}^2$ is 
a rank-3-bundle. Then 
$$D.c_2(X)>0 \mbox{ for all } D\in\overline{K(X)}.$$ 
\end{citethm}

Thus we are confirming the Wilson conjecture for hypersurfaces in ${\bb P}^2$-bundles over ${\bb P}^2$
with $c_3(X)\le -54$ and $\rho(X)=2$. 
The positivity of $c_2(X)$ is also proved by using the bundle information widely.

The cases $(n,m)=(3,1)$ and $(n,m)=(1,3)$ are dealt with in a different paper. Moreover, the case of a family of quartic elliptic curves may be also interesting.

Calabi-Yau manifolds with $\rho(X)=2$ are also considered in Mirror Symmetry as mirrors of two-parameter-families of Calabi-Yau manifolds, which are constructed 
to be able to choose a good one-parameter family.

This article grew out of the author's thesis at the University of Bayreuth. I would like to express my gratitude to Prof.\,Thomas Peternell for all his support and
valuable advice.
 
\section{Notation}

In this section we summarize the most important notations of this paper. $X$ will always denote a Calabi-Yau-threefold, while $Z$ is
a fourfold in every case. 
\begin{center}
\begin{tabular}{|l|l|}
\hline
$N^1(X)$& the ${\bb R}$-vector space of numerical classes of $Div(X)\otimes{\bb R}$\\
$K(X)$& the open part of the K\"ahler cone of $X$, i.e. the ample cone\\
$B(X)$& the effective cone in $N^1(X)$\\
$W(X)$& the hypersurface $\{D^3=0\}\subset N^1(X)$\\
${\cal O}_X(1)$ & the restriction ${\cal O}_Z(1)|X$, where ${\cal O}_Z(1)$ is the tautological bundle\\
& associated to $Z={\bb P}({\cal E})$\\
$K_Z$& the canonical divisor of $Z$\\
$c_i({\cal E})$& the $i$-th Chern class of a bundle ${\cal E}$\\
$c_i(M)$& the $i$-th Chern class of the tangent bundle of the complex\\
& manifold $M$\\
\hline
\end{tabular}
\end{center}

\section{The K\"ahler cone of Calabi-Yau threefolds}

If $X\stackrel{i}{\subset} Z$ is a Calabi-Yau-threefold in a $4$-fold $Z$, we want to know something about the relation of their K\"ahler cones. 
The easiest relation would be
$i^*K(Z)=K(X)$. This situation we will denote by $K(X)=K(Z)|X$. In \cite{bo89} Koll\'ar proved such a result, if $Z$ is Fano:

\begin{lemma}{\bf(Koll\'ar)} If $Z$ is a Fano-4-fold and $X\in|-K_Z|$ a Calabi-Yau-manifold, then 
$i^*:N^1(Z)\seq N^1(X)$
is an isomorphism and $K(X)=K(Z)|X$.
\end{lemma}

One could hope that this result holds also, if $-K_Z$ is only big and nef instead of ample. This is wrong, as the following example shows.

\begin{ex}{($K(X)\not=K(Z)|X$)}\label{ggbsp} Let ${\cal E}:={\cal O}\oplus{\cal O}(1)\oplus{\cal O}(2)\seq{\mathbb P}^2$ and $Z:={\mathbb P}({\cal E})$ 
with tautological bundle ${\cal O}_Z(1)$ and natural projection $p:Z\seq{\mathbb P}^2$. Then
$-K_Z={\cal O}_Z(3)$ is nef, globally generated and since 
$$-K_Z^4=81(c_1^2({\cal E})-c_2({\cal E}))=567>0,$$
$-K_Z$ is also big. The morphism
$$\Phi_{|-K_Z|}:Z\seq Z'$$
contracts $G:={\mathbb P}({\cal O})\cong{\mathbb P}^2$ to a point $z$, where ${\cal O}$ is the first summand of ${\cal E}$, and is an isomorphism outside $G$.
Since $-K_Z=\Phi^*{\cal L}$ for a globally generated line bundle ${\cal L}$ on $Z'$, we can choose a smooth element $L\in|{\cal L}|$, such that
$$X:=\phi^{-1}(L)$$
is smooth and $z\notin L$. 

Now, if $C\subset Z$ is a curve with ${\cal O}_Z(1).C=0$, then $C\subset G$. But since $X\cap G=\emptyset$, it follows $C\cap X=\emptyset$. This proves that
${\cal O}_X(1):={\cal O}_Z(1)|X$ is ample and hence
$$K(X)\not=K(Z)|X.$$

Moreover, it is not hard to proof that
$$\overline{K(X)}={\mathbb R}^+({\cal O}_X(1)-\pi^*h)+{\mathbb R}^+\pi^*h,$$
if we denote $\pi:=p|X$.
\end{ex}

One of the major steps to the rationality result, Theorem \ref{genkollar}, states, that this example already shows the structure of all counterexamples in the case $(2,2)$.
But for the moment we are not going to specialize the situation. Instead, we want to see what happens, if $-K_Z$ is not nef. Indeed, the statement of the 
Koll\'ar Lemma holds for the face of $K(Z)$ 'closer' to $-K_Z$:

\begin{lemma}\label{canside}Let $Z$ be a 4-fold with $\rho(Z)=2$ and assume $-K_Z$ is not nef. Moreover, let $X\in|-K_Z|$ be a Calabi-Yau threefold with $\rho(X)=2$ and $H\in Pic(Z)$ ample. If $r\in\bb R$ and $r>0$, then
$$H+r(-K_Z)\in\partial K(Z)\iff H+r(-K_Z)|X\in\partial K(X).$$
\end{lemma}

\begin{proof}Let $E_1,E_2$ be cone generators of $\overline{K(Z)}$- i.e. generators of $N^1(Z)$ with
$E_1,E_2\in\partial K(Z)$-, such that
$$-K_Z=E_1-E_2.$$
Furthermore let $c\in{\mathbb R}$ such that
$$D:=E_1-cE_2$$
is nef restricted to $X$. We get
$$-K_Z=D+(c-1)E_2.$$
Hence $c\ge 1$ would imply that $-K_Z|X$ is nef, therefore $-K_Z$ would be nef, what contradicts to the assumptions.

So we know $c<1$. Now let $C\subset Z$ be a curve. If $C\subset X$, then by construction, $D.C\ge 0$. 

If $C\not\subset X$, then
$(-K_Z).C\ge 0$, since $X\in|-K_Z|$. Therefore
$$D.C=-K_Z.C+(1-c)E_2.C\ge 0.$$
Hence $D$ is nef. This means that $D=E_1$. 

If now $H=h_1E_1+h_2E_2$ is ample and $H+r(-K_Z)|X\in\partial K(X)$ for $r>0$, these considerations imply that
$$H+r(-K_Z)=(h_1+r)E_1+(h_2-r)E_2=(h_1+h_2)E_1$$
and hence $H+r(-K_Z)\in\partial K(Z)$.
\end{proof}

Now the question is, if we can find ${\bb Q}$-Divisors on $\partial K(X)$ and, moreover, if we can express every element of $K(X)$ as a
finite sum of ${\bb R}$-linear combinations of ${\bb Q}$-divisors in $\partial K(X)$. The last property is called the rational generatedness (or short: rationality)
of the K\"ahler cone and is only dependent on the lattice $Pic(X)\subset N^1(X)$, but not on a chosen basis of $Pic(X)\otimes{\bb Q}$. 
 
First, we give a situation, in which we can use Lemma \ref{canside} to conclude 
the existence of a rational ray in $K(X)$.

\begin{lemma}\label{notnef}Let $Z$ be a 4-fold, $\rho(Z)=2$, $h^0(-K_Z)>1$ and $X\in|-K_Z|$ a Calabi-Yau threefold with $\rho(X)=2$. Then $\partial K(Z)$ contains a rational divisor of the form
$-K_Z+sH$ with $H\in Pic(Z)$ being ample and $0\le s\in{\bb Q}$. 
\end{lemma}

\begin{proof}
We divide the proof in two cases: $-K_Z$ nef and $-K_Z$ not nef.

Let us at first assume that $-K_Z$ is nef. 
If $-K_Z\in\partial K(Z)$, then we choose $s=0$. Hence we assume that $-K_Z$ is ample. 
Now we use the cone theorem for Fano manifolds, which states the rationality of 
$\partial K(Z)$.

If $-K_Z$ is not nef, we apply Lemma \ref{canside} and see that 
$$H+r(-K_Z)\in\partial K(Z)\iff H+r(-K_Z)|X \in\partial K(X),$$
if $r>0$ and $H\in Pic(Z)$ ample.
Since we assumed $h^0(-K_Z)>1$, we conclude by $X\in|-K_Z|$, that $-K_Z|X$ is effective. Now we want to use
the log-rationality theorem (\cite[Thm 4-1-1]{kmm}), which states that
$$\sup\{r\in\bb R|H+r(K_X+\Delta)\in K(X)\}\in\bb Q,$$
if $H\in Pic(X)$ is ample, $K_X+\Delta$ not nef and $\Delta$ an effective ${\bb Q}$-Divisor, 
such that $(X,\Delta)$ has only weak log-terminal singularities. The latter property is
guaranteed by choosing $\varepsilon\Delta$ instead of $\Delta$ for $0<\varepsilon\ll 1$.
Note that $K_X=0$.

Because $-K_Z|X$ is not nef, we can apply the log-rationality theorem on $X$
for arbitrary ample $H$ and $\Delta:=\varepsilon(-K_Z|X)$ for $0<\varepsilon\ll 1$ and $\varepsilon\in\bb Q$
and get, that 
$\partial K(X)$ contains a rational ray, by Lemma \ref{canside} coming from
a rational ray of $\partial K(Z)$. We set $s:=\frac 1r$.
\end{proof}

Following a notation of P.M.H. Wilson we denote 
$$W(X):=\{D\in N^1(X)|D^3=0\}.$$
The effective cone in $N^1(X)$ shall by named $B(X)$.
 
By a Lemma of Wilson in \cite{wi} we get to the following statement:

\begin{thm}\label{kwb}Let $X$ be a Calabi-Yau threefold with $\rho(X)=2$. If $D\in\partial K(X)$ and
$L\in K(X)$, then the statements
\begin{enumerate}
\item $D\in\partial K(X)\setminus W(X)$,
\item $D-\epsilon L\in B(X)\setminus\overline{K(X)}$ for all $0<\epsilon\ll 1$,
\item $rD\in\partial K(X)\cap Pic(X)$ for some $r\in{\bb R}$
\end{enumerate}
satisfy
$$(i)\imp (ii)\imp (iii).$$
\end{thm}

We give a short proof and point out, which part was proved by Wilson.

\begin{proof}
Let $D\in\partial K(X)\setminus W(X).$ Since $D^3>0$, also 
$$D^2.(D+H')>0 \mbox{ and }D.(D+H')^2>0,$$
if $H'$ is chosen to be ample. If we set $H:=D+H'$, then $H$ is ample and satisfies
$$D^2.H>0 \mbox{ and }D.H^2>0.$$
Hence, in a small neighbourhood of $D$ for every $D^\prime\in\ {\bb Q}\otimes Pic(X)$ holds
$${D^\prime}^2.H>0 \mbox{ and }D^\prime.H^2>0$$
Now, by a lemma of Wilson \cite[Key Lemma]{wi} $mD^\prime$ is effective for $m\gg 0$, if $D'$ is not nef. 
This proves the first implication.

To prove the second one, we use Kawamata's Rationality Theorem. Let $D':=D-\epsilon L$ for an
$\epsilon\ll 1$. Then $(X,{\eta}D^\prime)$ has weak log-terminal singularities for
$0<{\eta}\ll 1$. If $H$ is an ample divisor, the rationality theorem implies the existence of 
$0<q\in{\bb Q}$, such that
$$D'':=H+qD'\in\partial K(X).$$
Hence, there are $m\in{\mathbb Z}$ such that $mD''\in Pic(X)\cap \partial K(X)$. 
Now, if ${D,\tilde D}\subset\partial K(X)$ generate $N^1(X)$, consider the induced scalar product
$$(aD+\tilde a\tilde D,bD+\tilde b\tilde D):=ab+\tilde a\tilde b.$$
By choosing $\epsilon$ small enough, we may assume $(D',D)>0$.
Since $D''\in\partial K(X)$ and
$$(D'',D)=(H,D)+q(D',D)>0,$$
we see, that there is a $r>0$ with $mD''=rD$.
\end{proof}

So for proving rationality of the K\"ahler cone, it is enough to look at the cases $D^3=0$ for a $D\in\partial K(X)$. This will play a crucial role in the
following. Since also $\pi^*h^3=0$ and $\pi^*h\in\partial K(X)$, if we additionally assume a fibration $\pi:X\seq S$,
we see that we can reduce then to the case $\partial K(X)\subset W(X)$. 
This in particular applies to our $X$ in ${\bb P}^2$-bundles over
${\bb P}^2$.

\begin{thm}\label{bquer=kquer}Let $X$ be a Calabi-Yau threefold with $\rho(X)=2$. Then
$$\overline{B(X)}=\overline{K(X)}\iff\partial K(X)\subset W(X).$$
\end{thm}

\begin{proof}In view toward Lemma \ref{kwb}, only '$\Leftarrow$' has to be shown. 
Therefore let $D,{D'}$ be cone generators of $\overline{K(X)}$ and $D^3={D'}^3=0$. 
We discriminate between two cases:
\begin{enumerate}
\item $D^2>0,{D'}^2>0$ and
\item $D^2>0, {D'}^2\equiv 0$. 
\end{enumerate}
In case i) we know $D^2{D'}>0$ and $D{D'}^2>0$, because else $D^2\equiv 0$ or ${D'}^2\equiv 0$, since
$D^3={D'}^3=0$. 
If we write 
$$E:=\kappa D+\lambda D',$$ 
with $\kappa,\lambda\in{\bb R}$, for the numerical class of an effective divisor,
take a look at 
$$E.D^2\ge 0 \mbox{ and }E.{D'}^2\ge 0.$$ 
These inequations turn to
$$\lambda D^2 {D'}\ge 0 \mbox{ and }\kappa {D'}^2D\ge 0.$$ 
Hence $\kappa,\lambda\ge 0$, i.e. $B(X)\subset\overline{K(X)}$.

In case ii) we look similarly at 
$$E.D^2\ge 0 \mbox{ and }E.D.{D'}\ge 0.$$ 
The first inequation implies by $D^2{D'}>0$ that 
$\lambda\ge 0$ and the second by ${D'}^2D=0$ that $\kappa\ge 0$ holds. 
So again $B(X)\subset\overline{K(X)}$.
 
Since a sufficient big multiple of an ample divisor gets effective, we also have the inclusion $K(X)\subset B(X)$. 
Hence we conclude $\overline{B(X)}=\overline{K(X)}$. 
\end{proof}

\section{Calabi-Yau-threefolds in ${\mathbb P}^2$-bundles over ${\bb P}^2$}

In the following we are considering Calabi-Yau manifolds $X$ of the form
$X\subset{\bb P}({\cal E})=:Z$, $X\in|-K_{{\bb P}({\cal E})}|$, with ${\cal E}$ being a rank-3-bundle over 
${\bb P}^2$. Let us denote $p:{\bb P}({\cal E})\seq {\bb P}^2$ the bundle projection and
$\pi:X\seq{\bb P}^2$ the restriction of $p$ to $X$. The hyperplane class in ${\bb P}^2$ is denoted by $h$, a fibre
of $p$ by $F$. Then
$\gamma({\cal E}):=c_1^2({\cal E})-3c_2({\cal E})$ is invariant under ${\cal E}\mapsto{\cal E}\otimes L$, with $L$ being a line bundle
on ${\bb P}^2$. Thus we also write $\gamma(Z)$. Later we will see, that $\gamma$ is 
a topological invariant of $X$. The line bundle ${\cal O}_Z(1)|X$ shall be denoted by ${\cal O}_X(1)$.

The following sequences are the basics of all proofs:

\begin{eqnarray}
&0\seq T_{Z|{\bb P}^2}\seq T_Z \seq p^* T_{{\bb P}^2}\seq 0&\label{c1}\\
&0\seq {\cal O}_Z\seq p^*({\cal E}^\vee)\otimes{\cal O}_Z(1) \seq T_{Z|{\bb P}^2}\seq 0&\label{c2}\\
&0\seq T_X\seq T_Z|X \seq N_{X|Z}\seq 0&\label{c3}\\
&c_1({\cal O}_Z(1))^3-p^*c_1({\cal E}).c_1({\cal O}_Z(1))^2+p^*c_2({\cal E}).c_1({\cal O}_Z(1))=0&\label{c4}
\end{eqnarray}
if ${\cal F}\seq{\bb P}^2$ and $L\seq{\bb P}({\cal E})$ are line bundles. 

Sequence (\ref{c2}) is the relative Euler sequence and equation (\ref{c4}) defines the Chern classes of
${\cal E}$.

\subsection{Intersection Theory and Picard Number}

By using these sequences, we can calculate the intersection theory on $Z$ resp.\,$X$:

\begin{lemma}\label{intersection}
Let $X\subset{\bb P}({\cal E})=:Z$ and ${\cal E}\seq{\bb P}^2$ a rank-3-bundle, $X\in|-K_{{\bb P}({\cal E})}|$ smooth. Then

\begin{enumerate}

\item $c_1({\cal O}_Z(1))^4=c_1^2({\cal E})-c_2({\cal E})$

\item $c_1({\cal O}_Z(1))^3.p^*h=c_1({\cal E}).h$

\end{enumerate}

\end{lemma}

\begin{proof}
Multiplying (\ref{c4}) by $p^*h$ yields (ii).
Multiplying (\ref{c4}) by ${\cal O}_Z(1)$ and using (ii) yields (i).
\end{proof}

The intersection theory and the Chern classes on $X$ are described by the following Lemma. In particular, we prove that $\gamma$ is a topological invariant
of $X$.

\begin{lemma}\label{xprod}Let $X\subset{\bb P}({\cal E})=:Z, {\cal E}\seq{\bb P}^2$ a rank-3-bundle, $X\in|-K_{{\bb P}({\cal E})}|$ smooth. Then

\begin{enumerate}
\item $c_1({\cal O}_X(1))^3=\gamma+c_1^2({\cal E})+3c_1({\cal E}).h,$

\item $c_1({\cal O}_X(1)^2).\pi^*h=2c_1({\cal E}).h+3,$

\item $c_1({\cal O}_X(1)).F=3,$

\item $c_1({\cal O}_X(1)).c_2(X)= 36+12c_1({\cal E}).h+2\gamma(Z),$

\item $\pi^*h.c_2(X)=36,$

\item $c_3(X)=-6\gamma-162;$ in particular $\gamma\ge -27$, if $\rho(X)=2$.
\end{enumerate}
\end{lemma}

\begin{proof}Tedious calculations. The last part $\gamma\ge -27$, follows by $c_3(X)=2(\rho(X)-h^{1,2}(X))\le 4$.\end{proof}

Finally, we want to give a criterion for our assumption of $\rho(X)=2$ to be true.

\begin{thm}\label{rho=2}Let $X\subset{\bb P}({\cal E})=:Z$ a Calabi-Yau threefold with ${\cal E}\seq{\bb P}^2$ being a rank-3-bundle. If $-K_Z$
is big and nef, then
$$\rho(X)=2+h^2({\cal E}^\vee\otimes{\cal E}).$$
In particular, $\rho(X)=2$, if $-K_Z$ is big and nef and the generic splitting type $(e_1,e_2,e_3)$ 
(with $e_1\le e_2\le e_3$) of ${\cal E}$ is not $(0,0,3)$ (in the normalization $c_1({\cal E}).h\in\{1,2,3\}$).
\end{thm}

\begin{proof}We look at the two sequences
\begin{eqnarray}&0\seq N_{X|Z}^\vee\seq\Omega_Z|X\seq\Omega_X\seq 0&\label{?}\end{eqnarray}
and
\begin{eqnarray}&0\seq\Omega_Z\otimes K_Z\seq\Omega_Z\seq\Omega_Z|X\seq 0.&\label{??}\end{eqnarray}

Our first aim is to show $H^i(T_Z)=H^i({\cal E}^\vee\otimes{\cal E})$ for $i>1$. For this purpose we calculate
$R^ip_*(p^*({\cal E}^\vee)\otimes{\cal O}_Z(1))={\cal E}^\vee\otimes R^ip_*{\cal O}_Z(1)=0$ for $i>0$. Therefore the
Leray spectral sequence yields 
$$H^i(p^*({\cal E}^\vee)\otimes{\cal O}_Z(1))=H^i({\cal E}^\vee\otimes{\cal E}).$$
Sequence (\ref{c2}) shows, that $H^i(p^*({\cal E}^\vee)\otimes{\cal O}_Z(1))=H^i(T_{Z|{\bb P}^2})$ for $i>0$,
since by $R^ip_*{\cal O}_Z=0$ for $i>0$ we get $H^i({\cal O}_Z)=H^i({\cal O}_{{\bb P}^2})=0$ for $i>0$.
For applying (\ref{c1}), we verify by the projection formula
$R^ip_*p^*T_{{\bb P}^2}=0$ for $i>0$. 
Hence again the Leray spectral sequence implies
$$H^i(p^*T_{{\bb P}^2})=H^i(T_{{\bb P}^2})=0$$
for $i>0$. This implies with (\ref{c1}) now, that
$$H^i(T_{Z|{\bb P}^2})=H^i(T_Z)$$
for $i>1$. Therefore
$$H^i(T_Z)=H^i({\cal E}^\vee\otimes{\cal E})$$
for $i>1$.

So we see that
$$H^i(\Omega_Z\otimes K_Z)=H^{4-i}(T_Z)^\vee=H^{4-i}({\cal E}^\vee\otimes{\cal E})^\vee$$
for $i<3$. In particular,
$$H^1(\Omega_Z\otimes K_Z)=H^3({\cal E}^\vee\otimes{\cal E})^\vee=0.$$
Since $H^{2,1}(Z)=0$ and $N_{X|Z}^\vee=K_Z|X$ the cohomology sequences of (\ref{?}) and (\ref{??}) contain 
\begin{eqnarray}&0\seq H^1(\Omega_Z)\seq H^1(\Omega_Z|X)\seq H^2({\cal E}^\vee\otimes{\cal E})^\vee\seq 0&\label{!}
\end{eqnarray}
resp.\,
\begin{eqnarray}&0\seq H^1(K_Z|X)\seq H^1(\Omega_Z|X)\seq H^1(\Omega_X)\seq H^2(K_Z|X).&\label{!!}\end{eqnarray}

Since $-K_Z$ was assumed to be big and nef, $H^i(K_Z|X)=0$ holds for $i<3$ and by (\ref{!!}) and (\ref{!}) therefore
$$\rho(X)=h^1(\Omega_Z|X)=\rho(Z)+h^2({\cal E}^\vee\otimes{\cal E}).$$

To prove the second part, we show at first, that 
$\rho(X)=2$, if the generic splitting type $(e_1,e_2,e_3)$ of ${\cal E}$ (with
$e_1\le e_2\le e_3$) satisfies the condition $e_3-e_1<3$. This we see in this way: an element $0\not=t\in 
H^2({\cal E}^\vee\otimes{\cal E})$ induces via Serre duality (after choice of a basis of
$H^0({\cal E}^\vee\otimes{\cal E}\otimes{\cal O}(-3))^\vee$) a section $s\in H^0({\cal E}^\vee
\otimes{\cal E}\otimes{\cal O}(-3)), s\not=0$. This section cannot vanish on the general line $L$, 
hence induces a section $0\not= s\in H^0({\cal E}^\vee\otimes{\cal E}\otimes{\cal O}(-3)|L)$. This is possible, only if $e_3-e_1\ge 3$.

Since $-K_Z$ is big and nef it follows that $-K_Z|p^*L$ is nef, in particular 
$$0\le -K_Z.{\bb P}_L({\cal O}(e_1))=3e_1+3-e_1-e_2-e_3=2e_1+3-e_2-e_3$$ 
and therefore $e_3-e_1\ge 3$ only, if $e_2=e_1, e_3=e_1+3$.
\end{proof}

\begin{ex}{($\rho(X)\not=2$)}\label{rho=4} Have a look at
${\cal E}=2{\cal O}\oplus {\cal O}(3)$. Then ${\cal E}$ is globally generated and since $-K_Z={\cal O}_Z(3)$,
also $-K_Z$ is globally generated. 
This proves the existence of smooth Calabi-Yau threefolds in ${\bb P}({\cal E})$. Because 
$(-K_Z)^4=81(c_1^2({\cal E})-c_2({\cal E}))=729>0$ (see Lemma \ref{intersection}),
$-K_Z$ is big and nef, hence ${\cal E}$ satisfies the assumptions of Theorem \ref{rho=2}. Furthermore,
$${\cal E}^\vee\otimes{\cal E}=2{\cal O}(-3)\oplus 5{\cal O}\oplus 2{\cal O}(3)$$
and hence $h^2({\cal E}^\vee\otimes{\cal E})=2$, what means $\rho(X)=4$.
\end{ex}

\begin{ex}{($\rho(X)=2$)}
Theorem \ref{rho=2} shows also, that $\rho(X)=2$, if ${\cal E}$ is one of the following:
\begin{enumerate}
\item ${\cal O}\oplus{\cal O}(a)\oplus{\cal O}(b), 0\le a\le b\le 2, a+b\le 3$, 
\item $T_{{\bb P}^3}|{\bb P}^2$,
\item $T_{{\bb P}^2}\oplus{\cal O}$. 
\end{enumerate}
\end{ex}

\subsection{The generalized Koll\'ar Lemma}

We get the following generalization of Koll\'ar's Lemma:

\begin{thm}\label{genkollar}Let $X\subset{\mathbb P}({\cal E})=:Z$, with ${\cal E}\seq{\mathbb P}^2$ being a rank-3-bundle, and $\rho(X)=2$. Then
$$K(X)=K(Z)|X$$
holds, unless: 

$-K_Z$ is big and nef, but not ample, $-K_Z|X$ is ample and there is a surface $G\subset Z$ such that $X\cap G=\emptyset$ and
$$[\mu G]=9{\cal O}_Z(1)^2-(6c_1({\cal E}).h+9){\cal O}_Z(1).p^*h+(9c_2({\cal E})+3c_1({\cal E}).h+9-2c_1^2({\cal E}))F$$
for a certain $\mu>0$.
\end{thm}

This Theorem is a crucial step on the way to the proof of the rationality of $\partial K(X)$. Its proof will be divided into an ascending sequence of lemmata.
For simplicity, let us denote

\begin{tabular}{ll}$(A)$&${\cal E}\seq{\mathbb P}^2$ is a rank-3-bundle, $Z:={\mathbb P}({\cal E})$, $X\subset|-K_Z|$ is a\\
&Calabi-Yau-manifold with $\rho(X)=2$ and $-K_Z$ is big and nef,\\
&but not ample.
\end{tabular}

\begin{lemma}In situation $(A)$ holds: If $\Phi_{|-mK_Z|}:Z\seq Z'$ contracts only finitely many curves, 
then these are smooth and rational.\end{lemma}

\begin{proof}Let $C$ be an irreducible contracted curve. Since $-mK_Z=\Phi^*{\cal O}_{Z'}(1)$ we have $-mK_{Z'}={\cal O}_{Z'}(1)$ and  $Z'$ has only 
canonical singularities, in particular they are rational. Hence $R^1\Phi_*{\cal O}_{Z}=0$ and by
$$0\seq{\cal I}_C\seq{\cal O}_Z\seq{\cal O}_C\seq 0$$
we also get $R^1\Phi_*{\cal O}_C=0$. By the Leray spectral sequence we conclude $H^1({\cal O}_C)=0$. Hence $C$ is smooth and rational. 
\end{proof}

\begin{lemma}\label{nocurves}In situation $(A)$ holds: The exceptional locus of $\Phi_{|-mK_Z|}:Z\seq Z'$ contains a two-dimensional component. 
\end{lemma}

\begin{proof}We assume, $\Phi$ contracts only (rational) curves and $C\cong{\bb P}^1$ is such a curve.
By the adjunction formula and $K_Z.C=0$ we get
$c_1(N_{C|Z})=c_1(K_C)=-2$. Now we compute
$$\chi(N_{C|Z})=3(1-g(C))+c_1(N_{C|Z})=1>0$$
and therefore $C$ deforms in $Z$, hence there is a contracted surface, which contradicts our assumption.
\end{proof}

\begin{lemma}\label{generators}In situation $(A)$ holds: 
$$H^4(Z,{\bb Z})=<F,{\cal O}_Z(1).p^*h,{\cal O}_Z(1)^2>.$$\end{lemma}

\begin{proof}By the K\"unneth formula we already know that $b_4(Z)=3$.

To show that $v_1:=F,v_2:={\cal O}_Z(1).p^*h$ and $v_3:={\cal O}_Z(1)^2$ are a ${\bb Z}$-Basis for $H^4(Z,{\bb Z})$ reduces to showing that
the matrix $A=(a_{ij})=(v_i.v_j)$ is invertible over ${\bb Z}$. 
But by Lemma \ref{intersection} we know, that
$$A=\left( \begin{array}{ccc}
0 & 0 & 1\\
0 & 1 & c_1({\cal E}).h\\
1 & c_1({\cal E}).h & c_1^2({\cal E})-c_2({\cal E})
\end{array} \right)$$
so the Lemma is proven.
\end{proof}

\begin{lemma}\label{g9}In situation $(A)$ holds: If $\Phi_{|-mK_Z|}|X:X\seq X'$ is an isomorphism, 
then the exceptional locus of $\Phi$ is two-dimensional and is contracted to a finite number of points.
If $G$ is such a contracted surface, then
$$[\mu G]=9{\cal O}_Z(1)^2-(6c_1({\cal E}).h+9){\cal O}_Z(1).p^*h+(9c_2({\cal E})+3c_1({\cal E}).h+9-2c_1^2({\cal E}))F$$
for a $\mu>0$. \end{lemma}

\begin{proof}Let $E\seq V$ be the exceptional locus of $\Phi$. Since $-K_Z=\Phi^*{\cal O}(\frac 1m)$ as a  ${\bb Q}$-divisor
we have
$$mX=\Phi^*H$$
for an  ample $H\subset Z'$ and $m>0$. 
Since $H$ is ample, $H$ intersects with each positive-dimensional component of $V$. 
This implies, that $\Phi|X$ can be 
an isomorphism only, if $\dim V=0$.

By Lemma \ref{nocurves} the exceptional locus of $\Phi$ is at least two-dimensional.

Let us assume there is a divisor $D$, which is contracted to a point by $\Phi$.
Then $D=\kappa{\cal O}_Z(1)+\lambda\pi^*h$ and therefore
\begin{eqnarray}-K_Z.D.F&=&3\kappa\nonumber\\
-K_Z.D.{\cal O}_Z(1).\pi^*h&=&3\kappa c_1({\cal E}).h+\kappa(3-c_1({\cal E}).h)+3\lambda.\nonumber\end{eqnarray}
 
Because $D$ is contracted to a point,
$-K_Z.D\equiv 0$. Both equations imply now $\kappa=\lambda=0$. So there is only the case left, that 
$\Phi$ is small.

Now let $G\subset Z$ be an irreducible surface, which is contracted to a point by $\Phi$. 
By Lemma \ref{generators} we can write $G$ as 
$$G=\kappa{\cal O}_Z(1)^2+\eta{\cal O}_Z(1).p^*h+\lambda F.$$ 
Since $\Phi|X$ is assumed to be an isomorphism, this means, that
$-K_Z.G\equiv 0.$
This leads to equations $-K_Z.G.{\cal O}_Z(1)=0$ and $-K_Z.G.p^*h=0$, or in terms of $\kappa,\eta,\lambda$
\begin{eqnarray}(2c_1^2({\cal E})-3c_2({\cal E})+3c_1({\cal E}).h)\kappa + (2c_1({\cal E}).h+3)\eta + 3\lambda & = & 0\nonumber\\
(2c_1({\cal E}).h+3)\kappa + 3\eta & = & 0\nonumber\end{eqnarray}
Solving this system leads to
$$\mu G=9{\cal O}_Z(1)^2-(6c_1({\cal E}).h+9){\cal O}_Z(1).p^*h+[9c_2({\cal E})+3c_1({\cal E}).h+9-2c_1^2({\cal E})]F$$ 
for a $\mu\in{\mathbb R}$. Positivity of $\mu$ is a consequence of $G.F\ge 0$.
\end{proof}

Now we finally can prove Theorem \ref{genkollar}:

\begin{proof}[of Thm \ref{genkollar}] Let us assume $K(X)=K(Z)|X$ does not hold. Since $p^*h$ and $\pi^*h$ 
are nef and not ample, the properties 
$K(X)\not= K(Z)|X$ and $\rho(X)=2$ imply in view to Lemma \ref{canside} and the Koll\'ar Lemma that the divisor $-K_Z|X$ is ample and $-K_Z$ big and nef, but not ample. 
By 
$$0\seq (-K_Z)^{\otimes(m-1)}\seq (-K_Z)^{\otimes m}\seq (-K_Z)^{\otimes m}|X\seq 0$$
and $H^1((1-m)K_Z)=0$ (Kawamata-Viehweg), we see that
$$H^0(-mK_Z)\seq H^0(-mK_Z|X)$$
is surjective. Hence $\Phi_{|-mK_Z|_X|}=\Phi_{|-mK_Z|}|X$. Since $-K_Z|X$ is ample, this means that $\Phi_{|-mK_Z|}|X$ is an isomorphism.
Therefore we conclude by Lemma \ref{g9} the existence of a surface as described. 
\end{proof}

\subsection{The rationality of the K\"ahler cone}

As an immediate consequence of Lemma \ref{notnef} we get

\begin{lemma}\label{ratZ}Let ${\cal E}\seq{\bb P}^2$ be a rank-3-bundle and $Z:={\bb P}({\cal E})$. 
If $h^0(-K_Z)>1$ 
and there is a Calabi-Yau threefold $X\subset Z$ with $\rho(X)=2$,
then $\partial K({\bb P}({\cal E}))$ is rational.
\end{lemma}

Finally, we need a computational result on global sections:

\begin{lemma}\label{num}Let ${\cal E}\seq{\bb P}^2$ a rank-3-bundle, $X\subset{\bb P}({\cal E})=:Z$ a
Calabi-Yau threefold with $\rho(X)=2$. If ${\cal O}_X(1)$ is ample and $-K_Z$ is nef, then
\begin{enumerate}
\item $c_1({\cal E}).h\ge -1,$
\item $h^0({\cal O}_X(1)-\pi^*h)\ge \frac 13\gamma+\frac 16c_1^2({\cal E})+\frac 12c_1({\cal E}).h,$
\item if $\frac 13\gamma+\frac 16c_1^2({\cal E})+\frac 12c_1({\cal E}).h>0$, then $c_1({\cal E}).h\ge 1$.
\end{enumerate}
\end{lemma}

\begin{proof}
A standard computation yields 
$$0< {\cal O}_Z(1)^2.p^*h|X=3{\cal O}_Z(1)^3.\pi^*h+(3-c_1({\cal E}).h){\cal O}_Z(1)^2.F=2c_1({\cal E}).h+3,$$
what implies $c_1({\cal E}).h\ge -1$. Therefore the first part is proven.
By the nefness of $-K_Z={\cal O}_Z(3)+(3-c_1({\cal E}).h)p^*h$ we get furthermore, that
${\cal O}_Z(3)+4p^*h$ is nef. 

For a line $L$ in ${\bb P}^2$ we have the ideal sequence 
$$0\seq{\cal O}_X(1)-\pi^*h\seq{\cal O}_X(1)\seq{\cal O}_X(1)|\pi^*L\seq 0.$$
This implies
$$h^0({\cal O}_X(1)-\pi^*h)\ge h^0({\cal O}_X(1))-h^0({\cal O}_X(1)|\pi^*L).$$
The sequence
$$0\seq {\cal O}_Z(1)+K_Z|p^*L\seq{\cal O}_Z(1)|p^*L\seq{\cal O}_X(1)|\pi^*L$$
shows by $$H^0({\cal O}_Z(1)+K_Z|p^*L)=H^0({\cal O}_Z(-2)+(3-c_1({\cal E}).h)p^*h|p^*L)=0$$
and
$$H^1({\cal O}_Z(1)+K_Z|p^*L)=H^1(p_*({\cal O}_Z(-2)+(3-c_1({\cal E}).h)p^*h|p^*L))=0,$$ 
because $R^1p_*({\cal O}_Z(1)+K_Z|p^*L)=0$, that
$$H^0({\cal O}_Z(1)|p^*L)\cong H^0({\cal O}_X(1)|\pi^*L).$$

The lefthand side is easy to compute: since $L\cong{\bb P}^1$, ${\cal E}|L$ splits, hence we can write 
$${\cal E}|L={\cal O}(a)\oplus{\cal O}(b)\oplus{\cal O}(c) \mbox{ with }a\le b\le c.$$ 
Since ${\cal O}_Z(3)+4p^*h|p^*L$ is nef, ${\cal O}_Z(3)+6p^*h|p^*L$ is ample and therefore ${\cal O}_Z(1)+2p^*h|p^*L$ is ample, what is the definition of
${\cal E}\otimes{\cal O}(2)|L$ being ample. This implies $a\ge -1$. 

If $a\ge 0$, then $h^0({\cal O}_Z(1)|p^*L)=h^0({\cal E}|L)=
a+b+c+3=c_1({\cal E}).h+3$.  

If $a=-1, b\ge 0$, then $h^0({\cal O}_Z(1)|p^*L)=b+c+2=c_1({\cal E}).h+3$.

Finally, if $a=-1,b=-1,c>0$, then  $h^0({\cal O}_Z(1)|p^*L)=c+1=c_1({\cal E}).h+3$ again.

Because ${\cal O}_X(1)$ is ample, we compute
\begin{eqnarray}h^0({\cal O}_X(1))&=&\chi({\cal O}_X(1))\nonumber\\
&=&\frac 16{\cal O}_X(1)^3+\frac{1}{12}{\cal O}_X(1).c_2(X)\nonumber\\
&=&\frac 13\gamma+\frac 16c_1^2({\cal E})+\frac 32c_1({\cal E}).h+3\nonumber
\end{eqnarray}

Finally, we get

\begin{eqnarray*}h^0({\cal O}_X(1)-\pi^*h)&\ge& h^0({\cal O}_X(1))-h^0({\cal O}_X(1)|\pi^*L)\nonumber\\
&=&\frac 13\gamma+\frac 16c_1^2({\cal E})+\frac 32c_1({\cal E}).h+3-c_1({\cal E}).h-3\nonumber\\
&=&\frac 13\gamma+\frac 16c_1^2({\cal E})+\frac 12c_1({\cal E}).h\nonumber
\end{eqnarray*}

If the right hand side is positive, ${\cal O}_X(1)-\pi^*h$ is effective, and hence
$$0<{\cal O}_X(1).({\cal O}_X(1)-\pi^*h).\pi^*h=2c_1({\cal E}).h.$$
This proves the last part of the theorem.
\end{proof}

Now we are in a comfortable situation. If we assume $\partial K(X)$ to be not rational and $h^0(-K_Z)>1$, then by Lemma \ref{ratZ} we know that
$$K(X)\not=K(Z)|X.$$
Hence, by Theorem \ref{genkollar} we get, that 

\begin{tabular}{ll}
$(irr_+)$&$-K_Z$ is big and nef, not ample, $-K_Z|X$ is ample and\\ 
$(irr_G)$&there is a surface $G\subset Z$ with $X\cap G=\emptyset$ which is of\\
& cohomological class\\ 
&$[\mu G]=9{\cal O}_Z(1)^2-(6c_1({\cal E}).h+9){\cal O}_Z(1).p^*h$\\
&$+(9c_2({\cal E})+3c_1({\cal E}).h+9-2c_1^2({\cal E}))F$\\
&for a $\mu>0$.
\end{tabular}

By Lemmata \ref{kwb} and \ref{bquer=kquer} we furthermore can apply, that 

\begin{tabular}{ll}
$(irr_{eff})$&$\overline{B(X)}=\overline{K(X)}$ and $D^3=0$ if $D\in\partial K(X)$.
\end{tabular}

This is essentially the knowledge of the situation. For the sake of shortness, we denote

\begin{tabular}{ll}
$(irr)$&${\cal E}\seq{\mathbb P}^2$ is a rank-3-bundle, $Z:={\mathbb P}({\cal E})$, $X\subset|-K_Z|$ is a\\
&Calabi-Yau-manifold with $\rho(X)=2$, $h^0(-K_Z)>1$ and\\
&$\partial K(X)$ is not rational.
\end{tabular}

Furthermore, denote

\begin{tabular}{ll}
$(irr^0)$& $(irr)$ and ${\cal E}$ is normalized such that ${\cal O}_X(1)$ is ample,\\ 
& but ${\cal O}_X(1)-\pi^*h$ is not.
\end{tabular}

We are now going to apply the items above, which are implied by $(irr)$. 

\begin{lemma}\label{ratX20}$(irr)$ implies $-18<\gamma\le 1$.\end{lemma}

\begin{proof}Since $-K_Z$ is big and nef, we know
$$0<(-K_Z)^4=27\gamma+486,$$
what means $\gamma>-18$.
If we set $D:={\cal O}_Z(3)-kp^*h$ for a $k\in{\bb R}$ such that $D^3=0$ and $D\in\partial K(X)$ we compute
$$k=c_1({\cal E}).h+\frac 32-\sqrt{\frac 94-\gamma}.$$
Since there has to be a solution in ${\mathbb R}\setminus{\mathbb Q}$, we conclude $\gamma\le 1$.
\end{proof}

\begin{lemma}\label{tab}$(irr^0)$ implies $-1\le c_1({\cal E}).h\le 4$ and for every line $L\subset{\bb P}^2$ the type of ${\cal E}|L$ is 
contained in the following tabular:

\begin{center}
\begin{tabular}{|c|l|}
\hline
$c_1({\cal E})$ & ${\cal E}|L$ \\ \hline 
$-1$ & $(-1,-1,1)$, $(-1,0,0)$\\
$0$ & $(-1,0,1)$, $(0,0,0)$\\
$1$ & $(0,0,1)$\\
$2$ & $(0,0,2)$, $(0,1,1)$\\
$3$ & $(0,1,2)$, $(1,1,1)$\\
$4$ & $(1,1,2)$\\
\hline
\end{tabular}
\end{center}
\end{lemma}

\begin{proof}
By Lemma \ref{bquer=kquer} we know that 
$$H^0({\cal O}_X(1)-\pi^*h)=0.$$ 
The condition 
$$h^0(-K_Z)\ge 2$$
moreover implies then that
$$H^0({\cal O}_Z(1)-\pi^*h)=0.$$

Now we apply Lemma \ref{num} and get
$$h^0({\cal O}_X(1)-\pi^*h)\ge \frac 13\gamma +\frac 16c_1^2({\cal E})+\frac 12c_1({\cal E}).h.$$
Since we know by Lemma \ref{ratX20} that $-18<\gamma\le 1$, we conclude
$$h^0({\cal O}_X(1)-\pi^*h)>0,$$
if $c_1({\cal E}).h\ge 5$. By Lemma \ref{num} we know on the other hand
$c_1({\cal E}).h\ge -1.$
Hence only the cases
\begin{eqnarray}\label{c11}-1&\le c_1({\cal E}).h&\le 4\end{eqnarray}
remain.

Let $L\subset{\bb P}^2$ be a line. Then  
$${\cal E}|L={\cal O}(a)\oplus{\cal O}(b)\oplus{\cal O}(c)$$
with $a\le b\le c$. We call ${\cal E}|L$ 'of type $(a,b,c)$'.
Since $-K_Z$ is big and nef, we can compute
$$0\le -K_Z.{\bb P}({\cal O}(a))=3a+3-c_1({\cal E}).h$$
and therefore
\begin{eqnarray}a&\ge \lceil \frac 13c_1({\cal E}).h -1\rceil&\ge -1\label{a}\end{eqnarray}

Let us assume, $-K_Z.{\bb P}({\cal O}(b))=0$ holds. 
Then ${\bb P}({\cal O}(b))$ is contained in the exceptional locus of $\Phi_{|-mK_Z|}$. 
Since ${\bb P}({\cal O}(b))$ deforms in 
${\bb P}({\cal O}(a)\oplus{\cal O}(b))$, this implies that
$G:={\bb P}({\cal O}(a)\oplus{\cal O}(b))$
is a surface in the exceptional locus. Since $K(X)\not=K(Z)|X$, we can apply
Lemma \ref{g9} and get $-K_Z.G\equiv 0$.
We get a contradiction by computing
$$-K_Z.{\bb P}({\cal O}(a)\oplus{\cal O}(b)).p^*h=({\cal O}_{{\bb P}({\cal O}(a)\oplus{\cal O}(b))}(3)+(3-c_1({\cal E}).h)F).F=3,$$
where $F$ denotes the class of a fibre of
${\bb P}({\cal O}(a)\oplus{\cal O}(b))\seq L$.  

Hence $-K_Z.{\bb P}({\cal O}(b))>0$. Again we have 
$$-K_Z.{\bb P}({\cal O}(b))=3b+3-c_1({\cal E}).h.$$
Therefore
\begin{equation}\label{b}
\begin{split}c_1({\cal E})=0&\imp b\ge 0\\
c_1({\cal E})=3h&\imp b\ge 1
\end{split}
\end{equation}

The conditions (\ref{c11}),(\ref{a}),(\ref{b}) together with $a\le b\le c$ 
yield the tabulated cases.
\end{proof}

\begin{thm}\label{rat}Let ${\cal E}\seq{\mathbb P}^2$ be a rank-3-bundle, $Z:={\mathbb P}({\cal E})$, $X\subset|-K_Z|$ a
Calabi-Yau-manifold with $\rho(X)=2$ and $h^0(-K_Z)>1$. Then $\partial K(X)$ is rational.
\end{thm}

\begin{proof}
The statement is, that $(irr)$ and therefore $(irr^0)$ is impossible. So let us assume, $(irr^0)$ holds for some fixed ${\cal E}$, $X$ and $Z$. 
Then we can apply Lemma \ref{tab} and lead the cases $c_1({\cal E}).h=-1,...,4$ to a contradiction separately.

If $c_1({\cal E}).h\in\{1,4\}$, then by Lemma \ref{tab} we see, that
${\cal E}$ is uniform. This implies (e.g. by \cite{oss}) that
${\cal E}$ is 
\begin{itemize}
\item $T_{{\bb P}^2}(-1)\oplus{\cal O}$ or ${\cal O}\oplus{\cal O}\oplus{\cal O}(1)$ (for $c_1({\cal E}).h=1$), resp.\,
\item $T_{{\bb P}^2}\oplus{\cal O}(1)$ or ${\cal O}(1)\oplus{\cal O}(1)\oplus{\cal O}(2)$ (for $c_1({\cal E}).h=4$).
\end{itemize}
Both cases imply the same $Z$'s, hence only one of them could be true. But since
$-K_Z$ is ample for all those bundles, we get the rationality of $\partial K(X)$
by the Koll\'ar Lemma. This is a contradiction.

To exclude $c_1({\cal E})=-1$, we compute
$$\chi(m(-K_Z|X-4\pi^*h))=(\frac 92\gamma({\cal E})-9)m^3+(\frac 12\gamma({\cal E}) +6)m<0$$
for $m\gg 0$. 
Since ${\cal O}_X(3)=-K_Z|X-4\pi^*h$ is ample, we have
$$\chi(m(-K_Z|X-4\pi^*h))=h^0(m(-K_Z|X-4\pi^*h))>0,$$
what is a contradiction.

Similarly, if we assume $c_1({\cal E})=0$, then 
$-K_Z|X-3\pi^*h={\cal O}_X(3)$ is ample and we get again the condition
$$\chi(m(-K_Z|X-3\pi^*h))=\frac 92\gamma({\cal E}) m^3+(\frac 12\gamma({\cal E}) +9)m>0$$
for $m\gg 0$. This implies $\gamma({\cal E})\ge 0$. Since
$c_1({\cal E})=0$ and $\gamma({\cal E})\le 1$, it follows that 
$\gamma({\cal E})=0$.
But in this case $W=\{{\bb R}\pi^*h, {\bb R}{\cal O}_X(1), {\bb R}({\cal O}_X(1)-\pi^*h)\}$
is rational and by $(irr_{eff})$ we conclude that $\partial K(X)$ is rational.

In the case $c_1({\cal E})=2h$ we have that
$$-K_Z|X={\cal O}_X(3)+\pi^*h$$ 
and ${\cal O}_X(1)$ are ample, hence
$$\chi(m(-K_Z|X-\pi^*h))=(\frac 92\gamma({\cal E})+45)m^3+(\frac 12\gamma({\cal E}) +15)m>0$$
for $m\gg 0$, what implies 
$$\frac 92\gamma({\cal E}) +45\ge 0,$$
therefore $\gamma({\cal E})\ge -10$. 

Similarly to the cases $c_1=-1,0$ we compute
\begin{eqnarray*}h^0({\cal O}_Z(1))&=&h^0({\cal O}_X(1))\\
&=&\chi({\cal O}_X(1))\\
&=&\frac 16{\cal O}_X(1)^3+\frac{1}{12}{\cal O}_X(1).c_2(X)\\
&=&\frac 13\gamma({\cal E})+\frac{20}{3}\\
&\ge& 4
\end{eqnarray*}
and again we can conclude by $\overline{K(X)}=\overline{B(X)}$ 
and $h^0(-K_Z)>1$ like in the proof of Lemma \ref{tab} that
$$h^0({\cal O}_Z(1)-p^*h)=0.$$
In particular every $D\in|{\cal O}_Z(1)|$ is irreducible.

Now, according to $(irr_G)$ there is a surface $G\subset Z$, such that
$$-K_Z.G\equiv 0.$$

Since $-K_Z={\cal O}_Z(3)+p^*h$ and $G$ 
cannot contain any curve in a fibre (else
$0=-K_Z.C=p^*h.C$ for such a curve), 
for every curve $C\subset G$ holds
$${\cal O}_Z(1).C<0$$
and therefore $C\subset Bs({\cal O}_Z(1))$, where $Bs$ denotes the base locus
of a linear system.
If we now choose $D,D'\in|{\cal O}_Z(1)|$, then
$$G\subset Bs({\cal O}_Z(1))\subset D\cap D'.$$
 
Since $D\in|{\cal O}_Z(1)|$ is irreducible and $h^0({\cal O}_Z(1))>1$, 
it follows that $D\cap D'$ is twodimensional, hence 
$G$ is a component of $D\cap D'$.

According to Thm \ref{genkollar} the cohomological class of $G$ is
$$[\mu G]=9{\cal O}_Z(1)^2-3{\cal O}_Z(1).p^*h+(19-3\gamma({\cal E}))F$$
and therefore
$$9=\mu G.F\le \mu{\cal O}_Z(1)^2.F=\mu.$$
This implies $\mu=9$, because $G.F\in{\bb Z}$.
According to Lemma \ref{generators} the set 
$\{F,{\cal O}_Z(1).p^*h, {\cal O}_Z(1)^2\}$ is a ${\bb Z}$-basis
of $H^4(Z,{\bb Z})$, hence also holds
$$\mu|gcd(9,3,19-3\gamma({\cal E}))=1,$$
contradicting $\mu=9$. 

Let us turn our attention to the last case $c_1({\cal E})=3h$. Then
$-K_Z|X={\cal O}_X(3)$ is ample and by Lemma \ref{tab} the splitting type
of ${\cal E}$ can only be $(0,1,2)$ or $(1,1,1)$. 

Let $G$ be an exceptional surface as in Theorem \ref{genkollar} resp.\,Lemma \ref{g9}.
Since again
$$h^0({\cal O}_Z(1))=\chi({\cal O}_X(1))=\frac{\gamma}{3}+9\ge 4$$
by $\gamma>-18$, and the arguments, which show that $D\in|{\cal O}_Z(1)|$
is irreducible, apply here as well as in the case $c_1({\cal E})=2h$, we can
choose $D,D'\in|{\cal O}_Z(1)|$ which are irreducible and
intersect with $G$. But since ${\cal O}_Z(3)=\Phi^*H$ we also conclude, that
$$D=\Phi^*A, D'=\Phi^*A'$$
for ample divisors $A,A'$. Hence $D\cap G\not=\emptyset$ implies by
Lemma \ref{g9} that $G\subset D$, and therefore
$$G\subset D\cap D'.$$
As before we conclude that $\mu=9$ and $G$ cannot contain any curve in a 
fibre of $p$. This sums up to showing that $G$ is a section of $p$ and hence
given by an exact sequence
$$0\seq{\cal G}\seq{\cal E}\seq{\cal F}\seq 0$$
via $G={\bb P}({\cal F})$.
   
By the tangent sequence of $G$ and the property $-K_Z.G\equiv 0$ we compute ${\cal F}={\cal O}$.

But this implies that ${\cal O}_L$ is a quotient bundle if ${\cal E}|L$ for any line $L\subset{\bb P}^2$, what excludes the
case ${\cal E}|L=3{\cal O}(1)$. Hence ${\cal E}$ is uniform of splitting type $(0,1,2)$. This amounts to determining ${\cal E}$ as
\begin{itemize}
\item ${\cal O}\oplus{\cal O}(1)\oplus{\cal O}(2)$, or
\item $T_{{\bb P}^2}\oplus{\cal O}$, or
\item $T_{{\bb P}^2}(-1)\oplus{\cal O}(2)$, or
\item $S^2(T_{{\bb P}^2}(-1)).$
\end{itemize}

In the first case $\gamma=3$, contradicting $\gamma\le 1$. In the second and third case $\gamma=0$, hence $W(X)$ and therefore also $\partial K(X)$
are rational (arguments like above). In the fourth case, by using the Euler sequence and its derived sequences by taking $S^2$ of it, we compute
$\gamma=-9$. This we lead to a contradiction by the following arguments.

Let $D:={\cal O}_X(3)-k\pi^*h\in\partial K(X)$. Since $D^3=0$ we get like before
$$k=c_1({\cal E}).h+\frac 32-\sqrt{\frac 94 +9}=c_1({\cal E}).h+\frac 32-\frac 32\sqrt{5}<c_1({\cal E}).h-\frac 32=\frac 32.$$
Therefore for $D':={\cal O}_Z(2)-p^*h$ the restriction $D'|X$ is not nef, and because of $(irr_{eff})$ we conclude that it is also not effective. 
On the other hand we have 
$$H^0(D')=H^0(D'|X).$$
But
$$H^0(D')=H^0(S^2{\cal E}\otimes{\cal O}(-1))=H^0(S^2S^2T_{{\bb P}^2}\otimes{\cal O}(-5))$$
and we have an exact sequence
$$0\seq(\det{\cal H})^{\otimes 2}\seq S^2S^2{\cal H}\seq S^4{\cal H}\seq 0$$
for a rank-2-bundle ${\cal H}\seq{\bb P}^2$, where the projection is
given by
$$(t_it_j)(t_kt_l)\mapsto t_it_jt_kt_l.$$
This yields, 
tensored with ${\cal O}(-5)$ and applied to ${\cal H}=T_{{\bb P}^2}$,
$$0\seq{\cal O}(1)\seq S^2{\cal E}(-1)\seq S^4T_{{\bb P}^2}(-5)\seq 0,$$
what implies
$$h^0(D')=h^0({\cal O}(1))=3,$$
what contradicts to $(irr_{eff})$ by the previous calculation.

This finishes the proof.
\end{proof}

The only case for which we do not know the rationality of $\partial K(X)$, is $h^0(-K_Z)=1$. But we can give a bound $\Gamma$ for $\gamma$, such that
$h^0(-K_Z)>1$, whenever $\gamma\ge\Gamma$.

\begin{thm}\label{ratgamma}Let ${\cal E}\seq{\bb P}^2$ a rank-3-bundle, $X\subset{\bb P}({\cal E})=:Z$ a
Calabi-Yau threefold with  $\rho(X)=2$. If $\gamma(Z)\ge -18$, then $h^0(-K_Z)>1$, in particular $\partial K(X)$ is rational.
\end{thm}

\begin{proof}
The sequence
$$0\seq{\cal O}_Z\seq -K_Z\seq N_{X|Z}\seq 0$$
implies, that
\begin{eqnarray}h^0(-K_Z)>1&\iff& h^0(N_{X|Z})>0.\label{pos}\end{eqnarray}
The tangent sequence
$$0\seq T_X\seq T_Z|X\seq N_{X|Z}\seq 0$$
yields that
\begin{eqnarray}h^0(N_{X|Z})&\ge &h^0(T_Z|X)\label{k1}\end{eqnarray}
and
\begin{eqnarray}h^0(N_{X|Z})&\ge &h^1(T_X)-h^1(T_Z|X).\label{k2}\end{eqnarray}
Finally, the sequence
\begin{eqnarray}&0\seq T_Z\otimes K_Z\seq T_Z\seq T_Z|X\seq 0,&\label{tan}\end{eqnarray}
the relative tangent sequence  and the relative Euler sequence allow the statement
\begin{eqnarray}h^0(T_Z|X)&\ge& h^0(T_Z)-h^0(T_Z\otimes K_Z)\nonumber\\
&=&h^0(T_Z)-h^4(\Omega_Z)\nonumber\\
&=&h^0(T_Z)\nonumber\\
&\ge& h^0(T_{Z|{\bb P}^2})\nonumber\\
&\ge& h^0({\cal E}^\vee\otimes{\cal E})-1.\nonumber
\end{eqnarray}
If $h^0({\cal E}^\vee\otimes{\cal E})>1$, then by (\ref{k1}) and (\ref{pos}) also  $h^0(-K_Z)>1$ is proven.

So let us assume $h^0({\cal E}^\vee\otimes{\cal E})=1$. Now we have a closer look at (\ref{k2}).
The ideal sequence of $X\subset Z$, the relative tangent sequence of $Z\seq{\bb P}^2$ and the relative Euler sequence imply, that
\begin{eqnarray}h^1(T_Z|X)&\le& h^1(T_Z)+h^2(T_Z\otimes K_Z)\nonumber\\
&=&h^1(T_Z)+h^2(\Omega_Z)\nonumber\\
&=&h^1(T_Z)\nonumber\\
&\le& h^1(T_{Z|{\bb P}^2})+h^1(T_{{\bb P}^2})\nonumber\\
&\le& h^1({\cal E}^\vee\otimes{\cal E}).\nonumber
\end{eqnarray}

Since $h^0({\cal E}^\vee\otimes{\cal E})=1$, we get $h^2({\cal E}^\vee\otimes{\cal E})=0$ and 
$$h^1({\cal E}^\vee\otimes{\cal E})=1-\chi({\cal E}^\vee\otimes{\cal E})=-2\gamma-8.$$

Furthermore we know that
$$h^1(T_X)=2-\frac 12c_3(X)=3\gamma+83,$$
therefore
$$h^0(N_{X|Z})\ge 5\gamma +91>0,$$
since we assumed $\gamma\ge -18$.
\end{proof}

\begin{rem}\label{g>-18}If $-K_Z$ is nef, then 
$$0\le (-K_Z)^4=27\gamma+486,$$
i.e.
$$\gamma\ge -18.$$
Hence we can apply Theorem \ref{ratgamma} and get, that $\partial K(X)$ is rational. 
\end{rem}

\subsection{The positivity of $c_2(X)$}

In this last section we want to verify the positivity of $c_2(X)$ and therefore confirm the conjecture of Wilson for the Calabi-Yau threefolds considered in this section in the
case $h^0(-K_Z)>1$.

\begin{thm}\label{c2>0}Let $X\subset{\bb P}({\cal E})$ a Calabi-Yau threefold with $\rho(X)=2$, where ${\cal E}\seq{\bb P}^2$ is 
a rank-3-bundle. Then 
$$D.c_2(X)>0 \mbox{ for all } D\in\overline{K(X)}.$$ 
\end{thm}

\begin{proof}Again we choose ${\cal E}$ in such a way that ${\cal O}_X(1)$ is nef, but ${\cal O}_X(1)-p^*h$ 
is not.
Hence there is a $k\in\left[0;1\left[ \right.\right.$ such that $D:={\cal O}_X(1)-
k\pi^*h\in\partial K(X)$. We know $D^3\ge 0$. This is a condition, only if $\gamma\le 2$, and then this means
$$k\le \frac 13c_1({\cal E}).h+\frac 12-\frac 13\sqrt{\frac 94-\gamma},$$

In the case $\gamma\le 2$ we compute
$$D.c_2(X)=36+12c_1({\cal E}).h+2\gamma-36k\ge 18+2\gamma+12\sqrt{\frac 94-\gamma}.$$
Using $\gamma\ge -27$, we see $D.c_2> 0$ (the only zero of $18+2\gamma+12\sqrt{\frac 94-\gamma}$ is $\gamma=-54$).

If $\gamma\ge 3$ we distinguish between two cases: $-K_Z$ is nef and $-K_Z$ is not nef. 
Using the formulas of Lemma \ref{xprod}, we get: 
$$-K_Z|X.c_2(X)=6\gamma+216>0,$$ 
since $\gamma\ge -27$. By Lemma \ref{xprod} we already know
$$\pi^*h.c_2(X)=36>0.$$
If $-K_Z$ is not nef, then $-K_Z|X=N_{X|Z}$ is not nef. Hence $-K_Z|X+r\pi^*h\in\partial K(X)$ for a $r>0$ and hence
$D.c_2(X)>0$. 

So let us assume that $-K_Z$ is nef. We are able to restrict to $k>0$, since in the case $k=0$ we get $c_1({\cal E}).h\ge -1$ as in the proof
of Lemma \ref{num} and 
$$D.c_2=36+12c_1({\cal E}).h+2\gamma\ge 30>0.$$

If now $k>0$, also by $c_1({\cal E}).h\ge -1$ and $\gamma\ge 3$ we verify 
$$\frac 13\gamma+\frac 16c_1^2({\cal E})+\frac 12c_1({\cal E}).h>0$$
and hence get by Lemma \ref{num} that $c_1({\cal E}).h\ge 1$. Therefore
\begin{eqnarray}D.c_2&=&36+12c_1({\cal E}).h+2\gamma-36k\nonumber\\
& \ge & 48+2\gamma-36k\ge 12 + 2\gamma > 0.\nonumber
\end{eqnarray}

Since
$$\pi^*h.c_2(X)=36>0,$$
the proof is finished.
\end{proof}

\end{document}